\documentstyle[amstex,amscd,amssymb,theorem,12pt]{article}
\newcounter{assn}
\renewcommand{\theassn}{(\roman{assn})}
{\begin{list}{{\normalshape\theassn}\hskip 0.5em}%
{\usecounter{assn}\topsep=0ex \labelwidth=0em \leftmargin=0em
  \labelsep=0em \itemsep=0ex plus 0.2ex \parsep=\itemsep
  \partopsep=0ex}}%
{\end{list}}
\newtheorem{cla}{Claim}

\newtheorem{cor}[cla]{Corollary}
\newtheorem{lem}[cla]{Lemma}
\newtheorem{pro}[cla]{Proposition}
\newtheorem{thm}[cla]{Theorem}
\theorembodyfont{\rm}
\newtheorem{rem}[cla]{Remark}

\newtheorem{dfn}[cla]{Definition}

\newenvironment{pf}{\begin{ProofwCaption}{Proof}}{\end{ProofwCaption}}
\newenvironment{pf1}{\noindent {\em Proof of
proposition \ref{tnoc}.\/}}{\hfill \qedsymbol}
\newenvironment{ProofwCaption}[1]%
  {\addvspace\theorempreskipamount \noindent{\it #1.}\rm}%
  {\qed \par \addvspace\theorempostskipamount}
\newcommand{\qedsymbol}{\mbox{$\Box$}}
\begin{document}
\def\Box{\leavevmode\hbox{\vrule\vbox{\hsize=0.6em\hrule
         \hrule height0.4em depth0.2em width 0pt
         \hbox{\hskip0.6em plus 0pt minus 1fil}\hrule}\vrule}}
\def\qed{{\unskip\nobreak\hfil\penalty50
         \hskip2em\hbox{}\nobreak\hfil\Box
         \parfillskip=0pt \finalhyphendemerits=0 \par
         \medbreak\noindent\ignorespaces}}
\newcommand{\A}{{\cal A}}
\newcommand{\Alb}{{\operatorname{Alb}}}
\newcommand{\B}{{\cal B}}
\newcommand{\C}{{\Bbb C}}
\newcommand{\Q}{{\Bbb Q}}
\newcommand{\Coker}{{\operatorname{Coker}}}
\newcommand{\comp}{{\scriptstyle\circ}}
\newcommand{\D}{{\cal D}}
\newcommand{\diag}{{\operatorname{diag}}}
\newcommand{\Ext}{{\operatorname{Ext}}}
\newcommand{\Gr}{{\operatorname{Gr}}}
\newcommand{\Grass}{{\operatorname{Grass}}}
\newcommand{\h}{{\Bbb H}}
\newcommand{\Hom}{{\operatorname{Hom}}}
\newcommand{\id}{{\bf 1}}
\newcommand{\im}{{\operatorname{Im}}}
\newcommand{\Ker}{{\operatorname{Ker}}}
\newcommand{\into}{\hookrightarrow}
\newcommand{\cL}{{\cal L}}
\newcommand{\mwp}[2]{\pmatrix \id_{#1} & \psi\\ 0 & \id_{#2}
\endpmatrix}
\newcommand{\N}{{\cal N}}
\newcommand{\cO}{{\cal O}}
\newcommand{\PR}{{\Bbb P}}
\newcommand{\per}[2]{\pmatrix {#1} & \tau_2 & 1 & 0\\
     \tau_2 & {#2} & 0 & p \endpmatrix}
\newcommand{\sqper}[2]{\bmatrix {#1} & \tau_2 & 1 & 0\\
     \tau_2 & {#2} & 0 & p \endbmatrix}
\newcommand{\Pic}{{\operatorname{Pic}}}
\newcommand{\R}{{\Bbb R}}
\newcommand{\T}{T\langle 0 \rangle}
\newcommand{\X}{{X}}
\newcommand{\XX}[1]{Y_{{#1}}}
\newcommand{\Z}{{\Bbb Z}}
\newcommand{\ZmodN}{\Z_N}
\title{Degenerations of abelian surfaces and Hodge structures}
\author{K.~Hulek and J.~Spandaw}
\date{}
\maketitle
%
%
\section{Introduction}

The starting point of this note is twofold: In \cite{HKW} the first
author together with Kahn and Weintraub constructed and described a
torodial compactification of the moduli space of
$(1,p)$-polarized abelian surfaces with a (canonical) level structure.
Moreover, Mumford's construction was used to associate to each
boundary point a degenerate abelian surface. On the other hand
Carlson, Cattani and Kaplan gave in \cite{CCK} an interpretation of
torodial compactifications of moduli spaces of abelian varieties in
terms of mixed Hodge structures. Here we want to discuss a connection
between \cite{HKW} and \cite{CCK}. More precisely, we restrict
ourselves to corank 1 degenerations in the sense of \cite{HKW}, or
equivalently to type II degenerations, i.e. to cycles of elliptic
ruled surfaces. Our main result says that a degenerate abelian
surface associated to a boundary point is (almost) completely
determined by the boundary point (for a precise formulation see
theorem \ref{T}). The crucial ingredient in the proof is the Local
Invariant Cycle theorem which relates the variation of Hodge
structure (VHS) to the mixed Hodge structure (MHS) on the singular
surface.

In section \ref{section2} we collect some  basic facts about
semi-stable degenerations  of abelian surfaces, some of which are
well known. The MHS of a cycle of elliptic ruled surfaces is computed
in section \ref{section3}. Section \ref{section4} starts with a
review of corank 1 boundary points in moduli spaces of
$(1,p)$-polarized abelian surfaces with a level structure. Next the
VHS of a specific family associated to such a boundary point is
described. Using the Local Invariant Cycle theorem,
 this and the result from section \ref{section3} give theorem \ref{T}.

We would like to thank J. Steenbrink for helpful discussions.

Both authors are grateful to the DFG for financial support under
grant HU 337/2-4.

\section{Preliminaries on degenerations}\label{section2}

In this section we want to summarize basic facts on degenerations of
abelian surfaces some of which are well known (\cite{FM,P}).

\subsection{Cycles of elliptic surfaces} Let $\Delta=\{z\in\C:
|z|<1\}$ be the unit disk. We consider proper, flat families
$$
\begin{CD}
     \X\\
     @VV{p}V\\
     \Delta
\end{CD}
$$ with $X_t$ smooth abelian for $t\neq0$. We shall always assume the
total space $\X$ to be smooth and K\"ahler, and the components of the
singular fibre $X_0$ to be algebraic. We also assume $X$ to be
relatively minimal.

If $X_0$ has global normal crossings and no triple points then Persson
\cite[proposition~3.3.1]{P} or \cite[p.~11 and 17]{FM} has shown, that
$X_0$ is smooth abelian or a cycle of elliptic ruled surfaces, i.e.
$$ X_0=Y_1\cup\ldots\cup Y_N.
$$ The $Y_i$ are smooth ruled elliptic surfaces. $Y_i$ and $Y_{i+1}$
intersect transversally along a smooth curve which is a section of
both $Y_i$ and $Y_{i+1}$. In particular, the $Y_i$ are all ruled
surfaces over the same base curve $C$ with two disjoint sections. The
situation can be envisaged as in figure~1.
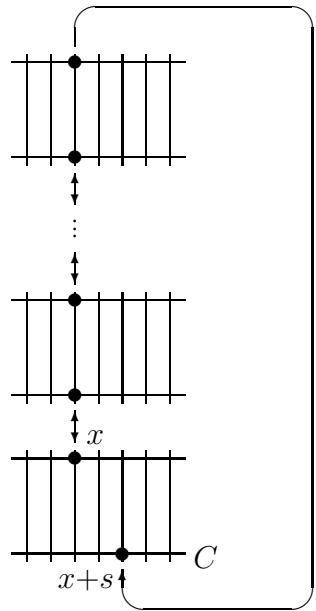
\begin{figure}
\newsavebox{\mybox}
\setlength{\unitlength}{3pt}
\begin{center}
  \begin{picture}(38,76)(10,5)
    \savebox{\mybox}(22,14){%
      \begin{picture}(22,14)
        \put(-2,12){\line(1,0){22}}
        \put(-2,0){\line(1,0){22}}
        \multiput(0,-1)(3,0){7}{\line(0,1){14}}
      \end{picture}}
    \put(12,62){\usebox{\mybox}}
    \put(18,74){\circle*{1.5}}
    \put(18,62){\circle*{1.5}}
    \put(18,56){\vector(0,1){4}}
    \put(18,56){\vector(0,-1){0}}
    \multiput(18,52)(0,1){3}{\circle*{0.1}}
    \put(18,46){\vector(0,1){4}}
    \put(18,46){\vector(0,-1){0}}
    \put(12,32){\usebox{\mybox}}
    \put(18,44){\circle*{1.5}}
    \put(18,32){\circle*{1.5}}
    \put(18,26){\vector(0,1){4}}
    \put(18,26){\vector(0,-1){0}}
    \put(12,12){\usebox{\mybox}}
    \put(18,24){\circle*{1.5}}
    \put(19.5,25.8){\makebox(0,0)[bl]{$x$}}
    \put(24,12){\circle*{1.5}}
    \put(23,10){\makebox(0,0)[tr]{$x{+}s$}}
    \put(33,11.5){\makebox(0,0)[l]{$C$}}
    \put(24,8){\vector(0,1){2}}
    \put(36,8){\oval(24,6)[b]}
    \put(48,8){\line(0,1){70}}
    \put(33,78){\oval(30,6)[t]}
    \put(18,78){\line(0,-1){2}}
    \savebox{\mybox}{}
  \end{picture}
\end{center}
\caption{A cycle of ruled surfaces over $C$ glued with a shift $s$}
\end{figure} We take $Y_1$ to be the bottom ruled surface and $Y_n$
the top ruled surface, and call $s$ the gluing parameter or shift of
$X_0$. A careful analysis of Persson's proof shows that, if one
replaces the hypothesis of global normal crossings by local normal
crossings, than the only additional case is the following: $X_0$ is
irreducible, its normalization is an elliptic ruled surface with two
disjoint sections and $X_0$ arises from this by gluing these two
sections with a shift
$s$. We shall always refer to this situation  (including the case
$N=1$) as a cycle of elliptic ruled surfaces.

The central fibre of such a degeneration is determined by the
following data:
\begin{enumerate}
\item the number $N$ of components
\item the base curve $C$
\item line bundles $\cL_i$ with
$\PR(\cO\oplus\cL_i)\cong \XX{i}$\footnotemark
\footnotetext{Strictly speaking for
$N=1$ we have to take the normalization.}
\item the gluing parameter $s$.
\end{enumerate}

We first want to show that these data are not independent of each
other. Before we can do this, we fix some more notation. First assume
that $N\ge 2$. We denote the intersection of
$\XX{i}$ and $\XX{i+1}$ by $C_i$. Furthermore we normalize the line
bundle $\cL_i$ in such a way that
$\cO_{\XX{i}}(C_{i-1})|_{C_{i-1}}=\cL_i$ and hence
$\cO_{\XX{i}}(C_{i})|_{C_{i}}=\cL_i^{-1}$. Under this assumption the
line bundle $\cL_i$ is uniquely determined. All indices have to be
read modulo $N$. The above discussion also makes sense in case $N=1$,
if we replace
$\XX{1}$ by its normalization. The next two propositions are special
cases of results proved by  Persson \cite{P}.

\begin{pro} $K_{\X}=\cO_{\X}$.
\end{pro}

\begin{pf} Since $K_{\X}|_{X_t}=\cO_{X_t}$ for $t\neq0$ it follows
that
$$
    K_X=\sum_{i=1}^N m_i \XX{i}
$$ where the $m_i$ are uniquely defined up to a common summand. After
possibly relabelling the components of $X_0$ we can assume $m_1$ to
be maximal. Moreover by adding multiples of a fibre of $p$ we may
assume that $m_1=0$. This already gives the result for $N=1$. Now
assume $N\geq2$.

Adjunction gives
$$
     (K_{\X}+\XX{i})|_{\XX{i}}=K_{\XX{i}}.
$$ Since
$$
     K_{\XX{i}}=-C_{i-1}-C_{i}+a_i f_{P_i}
$$ for a suitable ruling $f_{P_i}$ over a point $P_i\in C$ and
$$
     \XX{i}|_{\XX{i}} = -C_{i-1}-C_i
$$ this implies
$$
     (m_{i-1}-m_i) C_{i-1} +(m_{i+1}-m_i) C_i -a_i f_{P_i} =0.
$$ Since $m_1=0$ we find
$$
     m_0 C_0 +m_2 C_1 = a_1 f_{P_1}.
$$ Since, moreover, all $m_i\le 0$ this shows that $m_0=m_2=a_1=0$.
Continuing in this way we find $m_i=a_i=0$ for all $i$.
\end{pf}

\begin{pro}
$\operatorname{(i)}$ All line bundles $\cL_i$ are isomorphic.

\noindent $\operatorname{(ii)}$ $\deg \cL_i=0$ for all $i$.
\end{pro}

\begin{pf} (i) For the normal bundle of $C_i$ in $\X$ we have
\begin{align*}
     \N_{C_i/\X} &= \N_{C_i/\XX{i}} \oplus \N_{C_i/\XX{i+1}}\\
     &= \cL_i^{-1} \oplus \cL_{i+1}.
\end{align*} On the other hand, since $K_{\X}=\cO_{\X}$ and since
$C_i$ is an elliptic curve, we find again by adjunction that $\det
\N_{C_i/\X}=\cO_{C_i}$. This shows (i).

(ii) Assume that $\deg\cL_i=m\neq0$. Then $|(C_i)_{\XX{i}}^2|=m$  and
$(C_i)_{\XX{i+1}}^2=-(C_i)_{\XX{i}}^2$ for all $i$. But then a
topological argument (see \cite[p.~94]{P}) shows that
$$
     H_1(X_t,\Z) =\Z^3 \oplus \Z_m
$$ for general $t$. This contradicts the fact that $X_t$ is abelian.
\end{pf}

\subsection{Additional structures} We shall now consider further
structures on the family $p:X\to\Delta$. First of all we consider
degenerations of {\em polarized} abelian surfaces, i.e. we assume
that a line bundle $\cO_X(1)$ exists on $X$ such that
$\cO_{X_t}(1)=\cO(1)|_{X_t}$ for $t\neq0$ is a polarization on the
smooth abelian surface $X_t$, and that $\cO_{X_0}(1)$ is ample. We
are particulary interested in the case where $\cO_{X_t}(1)$
represents a polarization of type $(1,p)$ where
$p\geq3$ is a prime number.

We also want to assume that $p:X\to\Delta$ is a degeneration of
polarized abelian surfaces with a (canonical) level structure. For
the concept of (canonical) level structure on $(1,p)$-polarized
abelian surfaces see \cite[I.1]{HKW}. Before we give a formal
definition recall the following: Let
$Y=\PR(\cal O\oplus \cal L)$ be a $\PR^1$-bundle over an elliptic
curve $C$ with
$\deg\cal L=0$. Let $Y^0$ be the open part of $Y$ which is given by
removing the two sections defined by line bundles $\cal O$, resp.
$\cal L$. Then $Y^0$ is a
$\C^*$-bundle over the base curve
$C$. More precisely $Y^0$ carries the structure of a commutative
complex Lie group, and as such it is an extension of the form
$$ 1\longrightarrow\C^*\longrightarrow Y^0\longrightarrow
C\longrightarrow0.
$$ In fact $Y^0$ is a semi-abelian surface of rank 1.

\begin{dfn} Let $p:X\to\Delta$ be a degeneration of $(1,p)$-polarized
abelian surfaces. We say that this is a degeneration of
$(1,p)$-polarized abelian surfaces with a ({\em canonical}) {\em
level structure} if the following holds:

\noindent (i) There exists an open subset $X^0\subset X$ such that
$X^0\to\Delta$ is a family of abelian Lie groups with the following
property: $X_t^0=X_t$ for $t\neq0$ and $X_0^0$ is the smooth part of
a component of $X_0$ and as such carries the structure of a
semi-abelian surface.

\noindent (ii) There exists an action of $\Z_p\times\Z_p$ on $X$ over
$\Delta$ with the following properties: It leaves $\cO_X(1)$
invariant and defines a (canonical) level structure on $X_t$ for
$t\neq0$ (in particular it operates on $X_t$ by translation by
elements of order $p$). Moreover, the subgroup of $\Z_p\times\Z_p$
which stabilies $X_0^0$ acts on $X_0^0$ as a subgroup.
\end{dfn}

The two next results show that the presence of a polarization and a
level structure imposes strong conditions on the singular fibre $X_0$.

\begin{lem} Let $p:\X\to\Delta$ be a degeneration of abelian surfaces
with $(1,p)$-polarization and a (canonical) level-structure. Then
there are only two possibilities:

\noindent $\operatorname{(i)}$ The central fibre consists of one
component
$X_0$. If $\tilde{X_0}$ is its normalization, then
$\cO_{\tilde{X_0}}(1)=\cO_{\tilde{X_0}}(C_0 + p f_P)$ for a suitable
point $P$ in the base curve $C$.

\noindent $\operatorname{(ii)}$ The central fibre consists of $p$
components $\XX{i}$,
$i=1,\ldots,p$. In this case $\cO_{\XX{i}}(1)=\cO_{\XX{i}}(C_{i-1} +
f_{P_i})$ for suitable points $P_i$.
\end{lem}

\begin{pf} We have
$$
     \cO_{\XX{i}}(1) =\cO_{\XX{i}}(a_i C_{i-1} +b_i f_{P_i})
$$ for suitable integers $a_i$ and $b_i$. Since these line bundles
glue together to give a line bundle on $X_0$ it follows immediately
that all $b_i$ are equal. We denote this number by $b$. Since
$\cO_{X_0}(1)$ is ample, and since we are dealing with a degeneration
of abelian surfaces with a $(1,p)$-polarization it follows that
\begin{equation}\label{1}
     2b\sum_{i=1}^N a_i = 2p
\end{equation} with all $a_i>0$.

We now consider the action of $\Z_p\times\Z_p$ on the set
$\{\XX{1},\ldots,\XX{N}\}$. Let $G$ be the stabilizer of $\XX{1}$.
Since $p$ is a prime number, there are three possibilities:

(i) $G=\{1\}$. Then $N\ge p^2$ and this contradicts formula~(\ref{1}).

(ii) $G=\Z_p$. In this case $N\ge p$. It follows from
formula~(\ref{1}) that $N=p$ and $a_i=b=1$.

(iii) $G=\Z_p\times \Z_p$. The group $\Z_p\times \Z_p$ acts on $X_0$
as a group of automorphisms, hence it must leave its singular locus
invariant. Since $p$ is an odd prime number the group $G$ must leave
the curve $C_0$ invariant.

As a subgroup of $X_0^0$ the group $G$ acts by translation on $C_0$.
The multiplicative group
$\C^\ast$ contains no subgroup isomorphic to $\Z_p\times\Z_p$, hence
the group $G$ must contain at least a subgroup $\Z_p$ which acts
non-trivially on $C_0$. Since $\cO_X(1)$ is invariant under $G$ it
follows that the degree of
$\cO_X(1)$ restricted to $C_0$ must be divisible by $p$. I.e. $b$
must be divisible by $p$. By formula~(1) it follows that $b=p$,
$N=1$ and $a_1=1$.
\end{pf}

\begin{rem}  \normalshape It is easy to construct degenerations
$p:X\to\Delta$ of $(1,p)$-polarized abelian surfaces without a level
structure, such that the polarization on different components $Y_i$
is numerically different.
\end{rem}

In lemma~4 we have seen that the central fibre of a degeneration
$p:\X\to\Delta$ of $(1,p)$-polarized abelian surfaces with a
level-structure has either 1 or $p$ components. Recall that
$\XX{i}=\PR(\cO\oplus\cL_i)$ and that all $\cL_i$ are isomorphic and
of degree 0 (proposition~2). We shall denote this line bundle by
$\cL$. Moreover, assume that we have chosen an origin $O$ on the base
curve $C$. Then we can consider the shift $s$ as a point on $C$.

\begin{pro} Let $p:X\to\Delta$ be a degeneration of $(1,p)$-polarized
abelian surfaceswith a (canonical) level structure, and denote the
shift of $X_0$ by $s$. Then there are two possibilities:

\noindent $\operatorname{(i)}$ $N=1$ and $\cal L=\cO_C(ps-pO)$

\noindent $\operatorname{(ii)}$ $N=p$ and there exist a point $s'\in
C$ with $s=ps'$ such that $\cal L=\cO(s'-O)$.
\end{pro}

\begin{pf}  (i) Assume $N=1$. Then
$\cO_{\tilde{X_0}}(1)=\cO_{\tilde{X_0}}(C_0+pf_P)$ and hence
$$
     \cO_{\tilde{X_0}}(1)|_{C_0}=\cL\otimes \cO_C(pP),\qquad
     \cO_{\tilde{X_0}}(1)|_{C_1}=\cO_C(pP).
$$ Since $C_1$ and $C_0$ are glued with the shift $s$ a necessary and
sufficient condition for $\cO_{\tilde{X_0}}(1)$ to descend to a line
bundle on $X_0$ is
$$
     \cL\otimes\cO_C(pP)=\cO_C(pP)\otimes\cO_C(s-O)^{\otimes p}
$$ which gives the claim.

(ii) Assume $N=p$. Then
$\cO_{\XX{i}}(1)=\cO_{\XX{i}}(C_{i-1}+f_{P_i})$ and hence
$$
     \cO_{\XX{1}}(1)|_{C_1} =\cO_C(P_1),\qquad
     \cO_{\XX{2}}(1)|_{C_1}=\cL\otimes \cO_C(P_2).
$$ From gluing $\XX{1}$ and $\XX{2}$ we obtain
$$
     \cO_C(P_2)=\cL^{-1}\otimes \cO_C(P_1).
$$ Continuing in this way, we get
$$
     \cO_C(P_p)=\cL^{-(p-1)}\otimes \cO_C(P_1).
$$ Finally gluing $C^p$ and $C_0$ with a shift $s$ gives the condition
$$
     \cL \otimes \cO_C(P_1) =
          \cL^{-(p-1)}\otimes \cO_C(P_1)\otimes \cO_C(s-O)
$$ i.e.
$$
     \cL^p = \cO_C(s-O)
$$ as claimed.
\end{pf}

\begin{rem} \normalshape In \cite[part~II]{HKW} a number of explicit
examples of degenerations of
$(1,p)$-polarized abelian surfaces with a level structure were
constructed. It was shown that both types of degenerations which were
discussed above actually occur.  Moreover, all elliptic curves $C$
and all shifts $s$ can be realised.
\end{rem}

\section{The mixed Hodge structure on $H^1(X_0)$}\label{section3}

Let $X_0=\cup_{i\in\Z_N} Y_i$ be a cycle of elliptic ruled surfaces
as in the previous section.  In this section we calculate the MHS on
$H^1(X_0)$ and show that the base curve $C$ and the shift $s$  can be
recovered from it.

\subsection{The spectral sequence associated to $X_0$}

For technical reasons we want that $X_0$ has global normal crossings.
Therefore, we shall first  assume $N\geq2$. However, this is not an
essential hypothesis (see remark (\ref{R})).

The MHS on $H^q(X_0)$ can be computed as follows (see
\cite[p.103]{G}). Let $Y^0= \bigsqcup_{i\in\ZmodN} Y_i$ and
$Y^1=\bigsqcup_{i\in\ZmodN} C_i$, where $C_i=Y_i\cap Y_{i+1}$. The
maps $\alpha_i$ and $\beta_i$ are defined as the inclusions
$$
     \alpha_i:C_i\into Y_i\text{ and }\beta_i: C_i\into Y_{i+1}.
$$ Consider the double complex
\[
     A^{pq}=A^q(Y^p)
\] of global $C^{\infty}$ differential forms, where $d:A^{pq}\to
A^{p,q+1}$ is the exterior derivative and
$\delta: A^{0q}=\oplus_{i\in\ZmodN} A^q(Y_i) \to A^{1q}
=\oplus_{i\in\ZmodN} A^q(C_i)$ on $A^q(Y_i)$  is given by
$(-\beta_{i-1}^{\ast},\alpha_i^{\ast}):A^q(Y_i)\to A^q(C_{i-1})\oplus
A^q(C_i)$. These coboundary maps satisfy the relations
$d^2=0$, $d\delta=\delta d$ and $\delta^2=0$. There is a single
complex $(A^{\bullet}, D)$ associated to
$(A^{\bullet\bullet},d,\delta)$:
\[
     A^k=\oplus_{p+q=k} A^{pq}\quad\text{$D=(-1)^p d+\delta$ on
$A^{pq}$.}
\] We call
$
     \h^{k}:=H_D^{k}(A^{\bullet})
$ the hypercohomology of the double complex.

\begin{lem}\label{apq} The hypercohomology $\h^k$ is canonically
isomorphic to $H^k(X_0)$.
\end{lem}
\begin{pf} Let $i_p:Y^p\to X_0$ be the natural map. Let $\A^{pq}$  be
the sheaf on $X_0$ defined by
\[
     \A^{pq}(U)=A^q(i_p^{-1}U).
\] Set $\A^k:=\oplus_{p+q=k}\A^{pq}$ and let
$D$ be the sheafified version of the coboundary operator $D$ above.
It is shown in \cite[lemma~4.6]{GS}  that the complex
$(\A^{\bullet},D)$ is an acyclic resolution of the constant sheaf
$\C_{X_0}$, hence
$H^k(X_0,\C)=H^k(\Gamma(X_0,\A^{\bullet}))=H_D^k(A^{\bullet})=\h^k$.
\end{pf}

There exists a spectral sequence $E_r^{pq}$ with $E_0^{pq}=A^{pq}$
and $d_0=d$. The map $d_1$ is induced by $\delta$ and, more generally,
$d_r:E_r^{pq}\to E_r^{p+r,q-r+1}$ is induced by $D$. Notice that
$E_1^{pq}=H^q(Y^p)$. Since $A^{pq}=0$ for $p<0$ and $p>1$, the
spectral sequence degenerates at $E_2$, i.e.
$E_2^{pq}=E_{\infty}^{pq}$.

The weight filtration $W_{\bullet}$ on $\h$  is defined to be the
filtration induced by
$W_k(A^{\bullet})=\oplus_{q\le k} A^{\ast,q}$. By the theory of
spectral sequences we have
$E_{\infty}^{pq}=\Gr_{q}^W \h^{p+q}$. Since $A^{pq}=0$ for $p<0$ and
$p>1$, this boils down to a commutative diagram with exact rows
\[
\begin{CD}
     0 @>>> W_0 @>>> H^1(X_0)  @>>> \Gr_1^W @>>> 0\\
     @. @| @| @|\\
     0 @>>> E_{\infty}^{10} @>{i}>> \h^1 @>{\pi}>> E_{\infty}^{01}
@>>> 0.
\end{CD}
\] The map $i$ is induced by the inclusion $A^{10}\to A^{10}\oplus
A^{01}$; the map $\pi$ is induced by the projection $A^{10}\oplus
A^{01}\to A^{01}$.

The Hodge filtration on $\h^k$ is induced by
$F^p(A^k)=\oplus_{a+b=k} F^p(A^{a,b})$. In our case we have
\[
     F^1\h^1=\Ker(D|_{F^1(Y^0)})=\Ker(\delta: F^1(Y^0)\to F^1(Y^1)).
\]

\begin{lem}\label{wf} Let $X_0$ be a cycle of ruled surfaces over  a
smooth curve $C$. The weight filtration on $H^1(X_0)$ takes the form
$$
     0 \subset W_0 \subset W_1=H^1(X_0).
$$ Furthermore, $W_0$ is the unique 1-dimensional Hodge structure
$\T$ of type $(0,0)$ and $\Gr_1^W$ is canonically isomorphic to
$H^1(C)$ (as a Hodge structure).
\end{lem}

\begin{pf} Since $E_1^{pq}=H^q(Y^p)$ we have
$W_0=E_{\infty}^{10}=E_2^{10}=\Coker(\delta: H^0(Y^0)\to H^0(Y^1))$.
Now it is easy to see that
$$
\begin{CD}
     H^0(Y^0) @>{\delta}>> H^0(Y^1) @>>> \Z @>>> 0\\
     @. @| @AA{S}A\\
     @. \oplus H^0(X_i) @= \oplus \Z,
\end{CD}
$$ where $S$ is the summation map, is commutative and has an exact
top row, hence $W_0$ is indeed the trivial Hodge structure $\T$.

To show that
$\Gr_1^W=E_{\infty}^{01}=E_2^{01}=\Ker(\delta: H^1(Y^0)\to H^1(Y^1))$
is isomorphic to $H^1(C)$, we introduce the following maps: Let
$$
     \rho_i: Y_i\to C_i\text{ and }\sigma_i: Y_{i+1}\to C_i
$$ be the projections and let $\tau_i=\sigma_0 \comp
(\alpha_1\comp\sigma_1)
\comp(\alpha_2\comp\sigma_2)\comp\cdots\comp(\alpha_{i-1}\comp\sigma_{i-1}):
Y_i\to C_0=C$.  Then one easily shows that the following sequence is
exact
$$
\begin{CD}
     0 @>>> H^1(C) @>{\tau^{\ast}}>>  H^1(Y^0) @>{\delta}>> H^1(Y^1),
\end{CD}
$$ where $\tau^{\ast}:=(\tau_1^\ast,\ldots,\tau_N^\ast): H^1(C) \to
H^1(Y^0)=H^1(Y_1)\oplus\cdots\oplus H^1(Y_N)$.
\end{pf}

\subsection{The base curve and the extension
class}\label{subsection22}

The mixed Hodge structure $H^1(X_0)$ can be regarded as an extension
of the pure Hodge structures
$W_0$ and $\operatorname{Gr}^W_1$. We want to show how one can
recover the base curve $C$ and the shift $s$ from this extension. In
\cite{C} an object
$\operatorname{Ext}(\operatorname{Gr}^W_1,W_0)$ was introduced which
parametrizes all such extensions. In our case it is an abelian
variety, namely the base curve $C$.

\begin{lem}\label{lemma14} Let $W_{\bullet}$ be as above. Then there
exists a canonical isomorphism between
$\Ext(\Gr_1^W,W_0)$ and $\Alb(C)$.
\end{lem}

\begin{pf} For any Hodge structure $H$ we define $J^0 H =
H_{\C}/H_{\Z}+ F^0 H$. If $H_1$ and $H_2$ are two Hodge structures,
then we can apply this to $H=\Hom(H_1,H_2)$ and we get
$\Ext(H_1,H_2)=J^0 \Hom(H_1,H_2)$ (see \cite[prop.~2]{C}). Now if
$H_1=H^1(X)$, where $X$ is a smooth projective variety and $H_2=\T$,
then $\Hom(H_1,H_2)$ is the dual Hodge structure of $H^1(X)$ and has
weight $-1$. Its Hodge filtration has the form
$$
     0=F^1\subset F^0\subset F^{-1}=H^1(X)^\ast,
$$  where $F^0(H^1(X)^\ast)=(H^{0,1}(X))^\ast$. Hence
$H^1(X)^{\ast}/F^0(H^1(X)^{\ast})=H^{1,0}(X)^{\ast}$ and
$J^0(H^1(X)^{\ast})=H^{1,0}(X)^{\ast}/H_1(X,\Z)=\Alb(X)$.
\end{pf}

In \cite{C} we find the following algorithm to calculate the
extension class $e\in J^0\Hom(B,A)$ of an extension
$$
\begin{CD}
     0 @>>> A @>{i}>> H @>{\pi}>> B @>>> 0.
\end{CD}
$$
\begin{enumerate}
\item Choose an integral retraction $r:H \to A$, i.e. a map defined
over $\Z$ satisfying $r\comp i =\id_A$.
\item Define two Hodge filtrations on $A\oplus B$:
$$
    \tilde{F}_{\infty}^{\bullet}:=(r,\pi)(F^{\bullet}H)\text{ and }
    \tilde{F}_0^{\bullet}:=F^{\bullet}A \oplus F^{\bullet}B.
$$
\item Find a $\psi\in\Hom(B,A)_{\C}$ such that
$$
     \mwp{A}{B} \tilde{F}_0^{\bullet} =\tilde{F}_{\infty}^{\bullet}.
$$
\end{enumerate} Then $e=[\psi]\in J^0\Hom(B,A)$.

Let $X_0$ be a cycle of ruled surfaces over the elliptic curve
$C=V/\Lambda$ with shift $s\in \Alb(C)$. We now describe the ${\Bbb
Z}$-splitting of the exact sequence
$0\to W_0\to \h^1\to\Gr_1^W\to 0$, as required in the first step of
this algorithm, explicitely in terms of $C$ and $s$. Let
$\{\lambda_1,\lambda_2\}$ be a $\Z$-basis for
$\Lambda$ and let $x_i:V\to\R$ be the corresponding coordinate
functions ($i=1$, 2). Let $s_i=x_i(s)$, i.e.
\[
     s=s_1\lambda_1+s_2\lambda_2.
\] (The numbers $s_i$ are defined only up to integers.) Let $\phi\in
H^1(C)$. Since $\{dx_1,dx_2\}$ is a basis for
$H^1(C)$, we may represent $\phi$ by
$\xi_1 dx_1 + \xi_2 dx_2\in A^1(C)$, which we will also denote by
$\phi$. Consider
\[
     \tau^{\ast}\phi:=(\tau_1^{\ast},\ldots,\tau_N^{\ast})\phi\in
E_0^{01}.
\] Recall that $\tau^{\ast}\phi$ represents an element in
$E_{\infty}^{01}$. Set
\[
     f=(f_1,\ldots,f_N)\in E_0^{10},
\] where
\[
     f_i=\begin{cases} 0 & 1\le i<N\\
                     -  \xi_1 s_1 - \xi_2 s_2 & i=N
     \end{cases}
\] is a constant function on $C_i$.

\begin{lem} The pair $(\tau^{\ast}\phi,f)\in A^{01}\oplus A^{10}$
represents a class in $\h^1$.
\end{lem}

\begin{pf} First of all, $d(\tau^{\ast}\phi)=0$ since $d\phi=0$.
Furthermore,
$\delta(\tau^{\ast}\phi)=(0,\ldots,0,(t_{-s}-\id_C)^{\ast}\phi)$,
where $t_{-s}:x\mapsto x-s: C\to C$. Let $g=\xi_1 x_1+\xi_2
x_2\in\Hom_{\R}(V,\R)$. Then $\phi=dg$ and
$f_N=(t_{-s}-\id_C)^{\ast}g$.
\end{pf}

We write $(\tau^{\ast}\phi,f)=\sigma(\phi)$. By the lemma above, the
map $\sigma: H^1(C)\to E_0^{01}\oplus E_0^{10}$ induces a map
\[
     \bar{\sigma}: \Gr_1^W=H^1(C)\to \h^1.
\] Since $\pi:\h^1\to \Gr_1^W$ is induced by the projection
$A^{01}\oplus A^{10}\to A^{01}$, $\pi\comp\bar{\sigma}=\id_{\Gr_1^W}$,
i.e. $\bar{\sigma}$ is a splitting over $\R$. In fact, $\bar{\sigma}$
is defined over $\Z$:

\begin{lem}\label{fl} The map $\bar{\sigma}$ splits the sequence
\[
\begin{CD}
     0 @>>> W_0 @>{i}>> \h^1 @>{\pi}>> \Gr_1^W @>>> 0
\end{CD}
\] over $\Z$.
\end{lem}

Assuming this for the moment,  we show that the extension class is
identified with the shift.

\begin{cor}\label{cor17} The natural isomorphism
$\Ext(\Gr_1^W,W_0)\cong\Alb(C)$ identifies the extension class
$e=[\psi]$ with the shift $s$.
\end{cor}
\begin{pf} Let $r:\h^1\to W_0$ be the map such that
\[
     (r,\pi)=(i,\bar{\sigma})^{-1}:\h^1 \to W_0 \oplus \Gr_1^W.
\] As before we set
\[
    \tilde{F}_{\infty}^1=(r,\pi)F^1\h^1\text{ and }
   \tilde{F}_0^1=F^1(W_0)\oplus F^1 (\Gr_1^W)
\] in $W_0 \oplus \Gr_1^W$. The extension class is represented by a
map
$\psi: \Gr_1^W\to W_0$ such that
\[
     \begin{pmatrix} \id_{W_0} &- \psi\\
          0 & \id_{\Gr_1^W} \end{pmatrix}
          \tilde{F}_{\infty}^1=\tilde{F}_0^1.
\] Since $\tilde{F}_0^1\cap W_0=F^1(W_0)=0$ we get
\[
     r(\omega)-\psi(\pi(\omega))=0\quad\text{for all } \omega\in
F^1\h^1.
\]

Now assume that $C=\C/(\Z\tau+\Z)$ and let $\phi=\tau dx_1+dx_2$,
where $x_1$, $x_2$ are the coordinates dual to the basis
$\lambda_1=\tau$,
$\lambda_2=1$.    Then
$F^1H^1(C)=\C\phi$ (see \S4.2).  Recall that $(\tau^{\ast}\phi,0)\in
A^{01}\oplus A^{10}$ represents an element $\omega\in F^1\h^1$ such
that $\pi(\omega)=\phi$. Consider $\sigma(\phi)=(\tau^{\ast}\phi,f)$.
Since $i:W_0\to \h^1$ is induced by the inclusion $A^{10}\to
A^{01}\oplus A^{10}$,
\[
     \omega=i(-f)+\sigma(\phi),
\] hence
\begin{align*}
     \psi(\tau dx_1 + dx_2) &=\psi(\pi(\omega))\\
                            & = r(\omega)\\
                            &=-f
\end{align*} in $W_0$. Identifying $W_0$ with $\Z$  as in the proof
of lemma~\ref{wf} we find
\begin{align*}
     \psi(\tau dx_1+dx_2) &=- \sum f_i\\
          &=-f_N\\
          &=\tau s_1 + s_2\\
          &=\int_{0}^s \tau dx_1+dx_2.
\end{align*} Hence $\psi=\int_0^s$ on $F^1H^1(C)$, i.e.
$\Ext(\Gr_1^W,W_0)\cong\Alb(C)$ identifies the extension class
$e=[\psi]$ with the shift $s$.
\end{pf}

\begin{rem}\label{R}\normalshape The above result is also true for
$N=1$. In this case let $Y_0$ be the normalization of $X_0$. Let
$X_0'$ be the cycle of elliptic ruled surfaces consisting of two
copies of
$Y_0$ and glued with the same shift $s$ as $X_0$. There exists a
projection map $X_0'\to X_0$ given by contracting one of the
components of $X_0'$. By functoriality,  the isomorphism between the
vector spaces $H^1(X_0)$ and $H^1(X_0')$ induced by this map is in
fact an isomorphism of mixed Hodge structures. Hence we can work with
$X_0'$ instead of $X_0$ .
\end{rem}

\begin{rem}\normalshape Instead of cycles of elliptic ruled surfaces
we can also consider cycles of
$\PR^1$-bundles over an abelian variety $Z$. Again the Hodge
filtration takes on the form
$0=W_{-1}\subset W_0\subset W_1=H^1(X_0)$ where $W_0$ is the
1-dimensional Hodge structure of type
$(0,0)$. As before
$\operatorname{Gr}^W_1=H^1(Z)$ and hence
$\operatorname{Ext}(\operatorname{Gr}_1^W,W_0)=
\operatorname{Alb}(Z)$. The same proof as before also shows that the
extension class in
$\operatorname{Ext}(\operatorname{Gr}_1^W,W_0)$ again coincides with
the shift $s$.
\end{rem}

\subsection{Proof of lemma~\ref{fl}}

Let $c:[0,1]\to X_0$ be the loop described in figure~2
\cite[23.8]{Gre}.
\unitlength1cm
\begin{figure}\label{f2}
\begin{center}
\begin{picture}(12,5)
\thicklines
\put(1,1){\line(1,0){5}}
\put(6.4,1){\line(1,0){0.2}}
\put(6.9,1){\line(1,0){0.2}}
\put(7.4,1){\line(1,0){0.2}}
\put(8,1){\line(1,0){3}}
\put(1,5){\line(1,0){5}}
\put(6.4,5){\line(1,0){0.2}}
\put(6.9,5){\line(1,0){0.2}}
\put(7.4,5){\line(1,0){0.2}}
\put(8,5){\line(1,0){3}}
\put(1,1){\line(0,1){4}}
\put(3,1){\line(0,1){4}}
\put(5,1){\line(0,1){4}}
\put(9,1){\line(0,1){4}}
\put(11,1){\line(0,1){4}}
\thinlines
\put(1,2){\line(1,0){5}}
\put(6.4,2){\line(1,0){0.2}}
\put(6.9,2){\line(1,0){0.2}}
\put(7.4,2){\line(1,0){0.2}}
\put(8,2){\line(1,0){1}}
\put(9,2){\line(2,1){2.0}}
\put(9.2,2){\line(1,0){0.1}}
\put(9.6,2){\line(1,0){0.1}}
\put(10.0,2){\line(1,0){0.1}}
\put(10.4,2){\line(1,0){0.1}}
\put(10.8,2){\line(1,0){0.1}}
\put(1,2){\circle*{0.2}}
\put(3,2){\circle*{0.2}}
\put(5,2){\circle*{0.2}}
\put(9,2){\circle*{0.2}}
\put(11,2){\circle*{0.2}}
\put(11,3){\circle*{0.2}}
\put(0,2.2){\makebox(0.8,0.8)[rb]{$p_0$}}
\put(2,2.2){\makebox(0.8,0.8)[rb]{$p_1$}}
\put(4,2.2){\makebox(0.8,0.8)[rb]{$p_2$}}
\put(8,2.2){\makebox(0.8,0.8)[rb]{$p_{N-1}$}}
\put(11.2,3.2){\makebox(0.8,0.8)[lb]{$p_N=p_0$}}
\put(11.2,2.2){\makebox(0.8,0.8)[lb]{$p_N+s$}}
\put(1,0){\makebox(2,0.8)[t]{$Y_1$}}
\put(3,0){\makebox(2,0.8)[t]{$Y_2$}}
\put(9,0){\makebox(2,0.8)[t]{$Y_N$}}
\end{picture}
\caption{The loop $c$}
\end{center}
\end{figure}
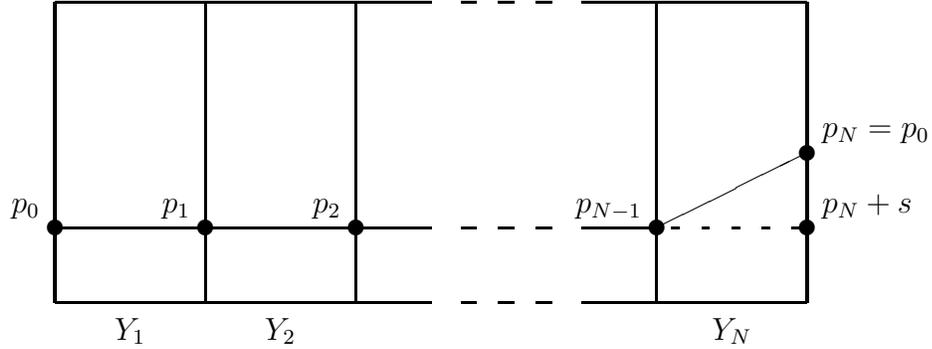 Clearly, $H_1(X_0,\Z)=H_1(C,\Z)\oplus \Z c$. Using the
Kronecker product (see \cite{Gre}) the loop $c$ determines a map
\[
     c: H^1(X_0,\Z)\to\Z.
\] Since $\bar{\sigma}$ splits the sequence over $\R$, all we have to
do is show that $c(\bar{\sigma}(dx_i))\in\Z$ for $i=1$, 2. In fact,
we will prove that $c\comp\bar{\sigma}=0$.

In order to understand how the loop $c$ operates on $\h^1$, we want
to identify $\h^1=H^1(\Gamma(X_0,\A^{\bullet}))$ with the simplicial
cohomology $H_{\nabla}^1(X_0)$. In \cite[p.95]{GM} we find a
canonical isomorphism
\[
     H^k(\Gamma(X_0,\B^{\bullet}))\to H_{\nabla}^k(X_0),
\] induced by integration, where
\[
     \B^k:=\Ker(\delta: \A^{0k}\to\A^{1k}).
\] One easily checks that the complex $(\B^{\bullet},d)$ is also an
acyclic resolution of the constant sheaf $\C_{X_0}$. The inclusions
$j_k:\B^k\into\A^{0k}\into\A^k$ commute with the coboundary
operators, hence the canonical isomorphism between
$H^k(\Gamma(X_0,\B^{\bullet}))$ and $H^k(\Gamma(X_0,\A^{\bullet}))$
is induced by $j_k$ (see \cite[5.24]{W}).

Let $\phi=\xi_1 dx_1+ \xi_2 dx_2$ and
$(\omega,f)=\bar{\sigma}(\phi)$, i.e.
\begin{align*}
     \omega
&=(\omega_1,\ldots,\omega_N)\qquad\omega_i=\tau_i^{\ast}\phi\\
     f&=(0,\ldots,0,f_N)\qquad f_N=-\xi_1 s_1- \xi_2 s_2
\end{align*} Let $g=(g_1,\ldots,g_N)\in\Gamma(X_0,\A^0)$ be such that
\[
     g_i|_{C_j}=\begin{cases} -f_N &(i,j)=(N,0)\\
                          0 &\text{otherwise}
                \end{cases}
\] Since $\delta g= -f$, we have that
$\omega'=\omega+dg\in \Gamma(X_0,\B^1)$.  Since
$(\omega',0)=(\omega,f)+(df,\delta f)$ we have
$j_1[\omega']=\bar{\sigma}(\phi)$ in $\h^1$.

Now we view the loop $c$ as an element of
$\Hom_{\Z}(H_{\nabla}^1(X_0,\Z),\Z)$. Using the defining properties
of $g_i$ and the fact that integrals of $\omega_i$ along the rulings
vanish,  we get
\begin{align*}
     c(\bar{\sigma}(\phi)) &=\sum_{i=1}^N \int_{p_{i-1}}^{p_i}
\omega_i'\\
                      &=\int_{p_{N-1}}^{p_N} \omega_N'\\
                      &=g_N(p_N)-g_N(p_{N-1})+\int_{p_N+s}^{p_N}
\omega_N\\
                      &=(s_1\xi_1 +s_2 \xi_2) - (s_1\xi_1 +s_2
\xi_2)\\
                      &=0.
\end{align*}
\hfill $\Box$

\section{Variation of Hodge structure}\label{section4}

In this section we want to compute the VHS associated to boundary
points of the moduli space of
$(1,p)$-polarized abelian surfaces with a (canonical) level structure
(cf. definition~\ref{D} and theorem~\ref{T}).

\subsection{$D$-polarized abelian varieties}

Our main reference is \cite{LB}. Fix a type $D$, i.e. an ordered
sequence $(d_1,\ldots,d_g)$  of positive integers satisfying
$d_i|d_{i+1}$ ($i=1,\ldots,g-1$). We will often write
$D=\diag(d_1,\ldots,d_g)$. Let $\h_g$ be the {\em Siegel space of
degree} $g$. To $\tau\in\h_g$ we associate an abelian variety
$X_{\tau,D}$ of dimension $g$ and a polarization $E_{\tau,D}$ of type
$D$ as follows. First,
\[
     X_{\tau,D}=\C^g/\Lambda_{\tau,D},
\] where
\[
     \Lambda_{\tau,D}=(\tau, D)\Z^{2g}.
\] Let $\lambda_i$ ($i=1,\ldots,2g$) be the $i$-th column of
$(\tau,D)$. Then we define $E_{\tau,D}$ to be the map
$\Lambda\times\Lambda\to\Z$ given by the matrix
\[
     E_{\tau,D}=\begin{pmatrix} 0 & D\\ -D & 0 \end{pmatrix}
\] with respect to the basis $\{\lambda_1,\ldots,\lambda_{2g}\}$ of
$\Lambda_{\tau,D}$. For this reason we will refer to this basis as
the symplectic basis. As an alternating map $E$ can be regarded as an
element in $H^2(X_{\tau,D},\Z)$. It is standard knowledge that (up to
sign) it is in fact the first  Chern class of an ample line bundle on
$X_{\tau,D}$, i.e. a polarization.

Every $D$-polarized abelian variety is isomorphic to
$(X_{\tau,D},E_{\tau,D})$ for some $\tau\in\h_g$. Furthermore, we can
construct a \lq\lq universal $D$-polarized abelian variety\rq\rq\
over $\h_g$ as follows:
\[
     \X_{D}=\Z^{2g}\backslash (\C^g\times\h_g),
\] where $\Z^{2g}$ operates on $\C^g\times\h_g$ by
\[
     l(v,\tau)=(v+(\tau,D)l,\tau)
\] for $l\in\Z^{2g}$ and $(v,\tau)\in\C^g\times\h_g$. The projection
$\C^g\times\h_g\to\h_g$ induces a projection
\[
     \pi:\X_D\to\h_g,
\] such that $\pi^{-1}(\tau)=X_{\tau,D}$.
$X_{D}$ is a complex manifold and carries a universal polarization
$\cL$ (see \cite[lemma~8.7.1]{LB}).

We are only interested in $D=(1,p)$, $D=(1)$ and $D=(p)$ and will
often omit $D$ if confusion seems unlikely.

\subsection{The polarized Hodge structure of
$(X_{\tau,D},E_{\tau,D})$}

We can describe the Hodge structure of the abelian variety
$X_{\tau,D}$ very explicitly. For any abelian variety $X=V/\Lambda$
we have
$$
\matrix
     F^1(H^1(X))=H^{1,0}(X) & \subset & H^1(X,\C)\\
     \| && \| \\
     V^{\ast}=\Hom_{\C}(V,\C) & \subset & \Hom_{\R}(V,\C).
\endmatrix
$$  Now consider $X=X_{\tau,D}$, let $\{e_1,\ldots,e_g\}$ be the
standard basis of $V=\C^g$ and let
$\{\lambda_1,\ldots,\lambda_{2g}\}$ be the symplectic basis of
$\Lambda=\Lambda_{\tau,D}$. Then
$\{\lambda_1^{\ast},\ldots,\lambda_{2g}^{\ast}\}$ is a $\Z$-basis of
$H^1(X,\Z)$ and $\{e_1^{\ast},\ldots,e_g^{\ast}\}$ is a $\C$-basis of
$F^1(H^1(X))$. Using the coordinates of $H^1(X,\C)=H^1(X,\Z)\otimes\C$
determined by $\{\lambda_1^{\ast},\ldots,\lambda_{2g}^{\ast}\}$,
$e_i^{\ast}$ is given by the $i$-th row of the period matrix
$(\tau,D)$.

A polarization $\omega$ on a compact complex manifold $X$ induces a
bilinear form $Q$ on $H^n(X)$:
$$
     Q(\phi,\psi)=(-1)^{n(n-1)/2} \int_{X}
\phi\wedge\psi\wedge\omega^{d-n},
$$ where $d=\dim X$. The pair $(H^n(X), Q)$ is a polarized Hodge
structure of weight~$n$ (see \cite[p.7]{G}).

We now calculate the polarization on $H^1(X_{\tau,D})$ induced by
$E_{\tau,D}$. Let $x_i:V\to \R$ ($i=1,\ldots,2g$) be the coordinates
with respect to the real basis $\{\lambda_1,\ldots,\lambda_{2g}\}$ of
$V$. Then $dx_i$ corresponds to $\lambda_i^{\ast}$ under the natural
isomorphism $H_{\hbox{de Rham}}^1(X) \cong \Hom(\Lambda,\Z)$.
Furthermore $E_{\tau,D}\in H^2(X,\Z)$ is represented by the 2-form
\[
     \omega=-\sum_{i=1}^g d_i dx_i\wedge dx_{i+g}
\] and
\[
     \int_X dx_1\wedge dx_{g+1} \wedge\cdots\wedge dx_{g}\wedge
dx_{2g} =1
\] (see \cite[lemmas 3.6.4 and~3.6.5]{LB}). It follows that $Q$ is
given by the matrix
$$
    (g-1)!
     \left(\begin{array}{cc} 0 & -\hat{D}\\ \hat{D} &
0\end{array}\right)
     \qquad  \hat{D}=\Big(\prod_{i=1}^{2g} d_i\Big) D^{-1}
$$ with respect to the basis
$\{\lambda_1^{\ast},\cdots,\lambda_{2g}^{\ast}\}$ of
$H^1(X,\Z)=\Hom(\Lambda,\Z)$.

\subsection{Boundary points of ${\cal A}^*(1,p)$}

We have to recall briefly some facts about compactifications of moduli
spaces of abelian surfaces. All relevant details can be found in
\cite{HKW}. By $\cal A(1,p)$ we denote the moduli space of
$(1,p)$-polarized abelian surfaces with a (canonical) level structure
($p\geq3$, prime). Recall that
$$
\cal A(1,p)=\h_2/\Gamma_{1,p}
$$ where
$$
\Gamma_{1,p}=\left\{g\in\operatorname{Sp}(4,\Z)\ ;\
g-\mbox{\boldmath$1$}\in
\begin{pmatrix}
\Z&\Z&\Z&p\Z\\ p\Z&p\Z&p\Z&p^2\Z\\
\Z&\Z&\Z&p\Z\\
\Z&\Z&\Z&p\Z
\end{pmatrix}\right\}
$$ acts on Siegel space $\h_2$ by
$$ g=\begin{pmatrix} A&B\\C&D\end{pmatrix}
:\tau\longmapsto(A\tau+B)(C\tau+D)^{-1}.
$$ In \cite{HKW} a torodial compactification ${\cal A}^*(1,p)$ of
$\cal A(1,p)$ was constructed. To compactify $\cal A(1,p)$ one has to
add (non-compact) boundary surfaces (corank 1 boundary points) and
boundary curves (corank 2 boundary points). Here we shall restrict
ourselves exclusively to the boundary surfaces. These surfaces are
indexed by the vertices of the Tits building of
$\Gamma_{1,p}$ which correspond to lines $l\subset\Q^4$. According to
\cite{HKW} there are two types of boundary surfaces, namely one {\em
central boundary surface} $D(l_0)$ and $p(p-1)/2$ {\em peripheral
boundary surfaces} $D(l_{(a,b)})$ where
$(a,b)\in(\Z_p\times\Z_p\setminus\{0\})/(\pm1)$. The group
$\operatorname{SL}(2,\Z_p)$ acts on ${\cal A}^*(1,p)$  and permutes
the peripheral boundary surfaces of ${\cal A}^*(1,p)$ transitively.
Therefore it is enough to consider one of them, namely $D(l_{(0,1)})$.

{}From now on let $l=l_0$ or $l_{(0,1)}$. The stabilizer $P(l)$ of $l$
in
$\Gamma_{1,p}$ is an extension of the form
$$ 1\longrightarrow P'(l)\longrightarrow P(l)\longrightarrow
P''(l)\longrightarrow 1
$$ where $P'(l)$ is a rank 1 lattice. The compactification procedure
requires that one first takes the partial quotient of $\h_2$ with
respect to $P'(l)$. For $l=l_0$
$$ P'(l)=\left\{\left(
\begin{array}{c|c}
  \begin{array}{ccc} && \\ &\mbox{\boldmath$1$}&\\ && \end{array} &
  \begin{array}{ccc} n&&0\\  && \\ 0&&0\end{array} \\ \hline
  \begin{array}{ccc} && \\ &0& \\ &&\end{array} &
  \begin{array}{ccc} && \\ &\mbox{\boldmath$1$}& \\ &&\end{array}
\end{array}\right)\ ;\ n\in\Z\right\}
$$ and for $l=l_{(0,1)}$:
$$ P'(l)=\left\{\left(
\begin{array}{c|c}
  \begin{array}{ccc} && \\ &\mbox{\boldmath$1$}&\\ && \end{array} &
  \begin{array}{ccc} 0&&0\\  && \\ 0&&np^2\end{array} \\ \hline
  \begin{array}{ccc} && \\ &0& \\ &&\end{array} &
  \begin{array}{ccc} && \\ &\mbox{\boldmath$1$}& \\ &&\end{array}
\end{array}\right)\ ;\ n\in\Z\right\}
$$ The partial quotient map $e(l):\h_2\to\h_2/P'(l)$ is then given by
$$ e(l_0):\h_2\longrightarrow\C^*\times\C\times\h_1,\
\begin{pmatrix}\tau_1&\tau_2\\ \tau_2&\tau_3
\end{pmatrix}\longmapsto(e^{2\pi i\tau_1},\tau_2,\tau_3)
$$ resp.
$$ e(l_{(0,1)}):\h_2\longrightarrow\h_1\times\C\times\C^*,\
\begin{pmatrix}\tau_1&\tau_2\\ \tau_2&\tau_3
\end{pmatrix}\longmapsto(\tau_1,\tau_2,e^{2\pi i\tau_3/p^2}).
$$ Here $\h_1$ denotes the usual upper half plane. This maps $\h_2$
to an interior neighbourhood of
$\{0\}\times\C\times\h_1$, resp. $\h_1\times\C\times\{0\}$, i.e. the
interior of the closure of the image, $X(l)=(\overline{e(\h_2)})^o$
is an open neighbourhood of
$\{0\}\times\C\times\h_1$, resp.
$\h_1\times\C\times\{0\}$. The partial compactification in the
direction of
$l$ then consists of adding the set $\{0\}\times\C\times\h_1$, resp.
$\h_1\times\C\times\{0\}$. There is a natural map
$X(l)\to{\cal A}^*(1,p)$ given by dividing out the extended action of
$P''(l_0)$ on $X(l)$, which maps $\{0\}\times\C\times\h_1$, resp.
$\h_1\times\C\times\{0\}$ to the boundary surface $D(l_0)$, resp.
$D(l_{(0,1)})$.

\subsection{Two 1-parameter families}\label{subsection34}

We consider two 1-parameter families of $(1,p)$-polarized abelian
surfaces which are closely related to the central, resp. peripheral
boundary components (see proof of theorem~\ref{T}). For $M$ a
positive integer we set
$$
\h_1(M):=\{\tau\in\h_1\ ;\ \operatorname{Im}\tau>M\}
$$ and
$$
\Delta^*(M):=\{t\in\C\ ;\ 0<|t|<e^{-2\pi iM}\}.
$$ First we fix a pair $(\tau_2,\tau_3)\in\C\times\h_1$. For
sufficiently large
$M$ we have a map
$$
\h_1(M)\longrightarrow\h_2,\ \tau\longmapsto\begin{pmatrix}
\tau&\tau_2\\\tau_2&\tau_3\end{pmatrix}.
$$ The first 1-parameter family which we want to consider is the pull
back of the universal family
$\pi:X_D\to\h_2$ to a family $\pi_1:X_D\to\h_1(M)$ via this map. For
the second family we fix a pair $(\tau_1,\tau_2)\in\h_1\times\C$ and
pull back the universal family
$\pi:X_D\to\h_2$ to a family $\pi_2:X_D\to\h_1(M)$ via the map
$$
\h_1(M)\longrightarrow\h_2,\
\tau\longmapsto\begin{pmatrix}\tau_1&\tau_2\\
\tau_2&p^2\tau\end{pmatrix}
$$ where $M$ is again chosen sufficiently large.

We denote the polarized abelian surface associated to $\tau\in\h_1(M)$
by $(X_{\tau},E_{\tau})$. Since in both cases $(X_{\tau},E_{\tau})$
depends only on $\tau\mod\Z$, the family $\pi_i$ is the pull back of
a family $p_i: \X_D\to\Delta^{\ast}(M)$ via the map
\[
\h_1(M)\longrightarrow\Delta^*(M),\ \tau\longmapsto t=e^{2\pi i\tau}.
\]

Consider the maps
\[
     \tilde{\phi}_1:\h_1(M)\to\Grass(2,4),\
\tau\mapsto\sqper{\tau}{\tau_3},
\] and
\[
     \tilde{\phi}_2:
          \h_1(M)\to\Grass(2,4),\ \tau\mapsto\sqper{\tau_1}{p^2\tau}.
\] There image is contained in Griffiths' period domain
$\D_2\subset\Grass(2,4)$. The monodromy $T_i$ is defined by
$\sideset{^t}{}{\tilde{\phi}}_i(\tau+1)=T_i
\sideset{^t}{}{\tilde{\phi}}_i(\tau)$. One easily checks that
$T_i=\left(\smallmatrix \id_2 & \eta_i\\  0 &
\id_2\endsmallmatrix\right)$, where $\eta_1=\left( \smallmatrix 1 &
0\\ 0 & 0 \endsmallmatrix \right)$ and $\eta_2=\left( \smallmatrix 0
& 0\\ 0 & p \endsmallmatrix \right)$ By construction $\tilde{\phi}_i$
descends to a map
$\phi_i: \Delta^\ast(M) \to \langle T_i \rangle \backslash \D_2$:
$$
\begin{CD}
     \h_1 @>{\tilde{\phi_i}}>> \D_2\\
     @V{t}VV @VVV\\
     \Delta^\ast(M) @>>{\phi_i}>\langle T_i \rangle \backslash \D_2
\end{CD}
$$ If we identify $H^1(X,\C)$ with $\C^4$ using the coordinates
induced by the basis $\{\lambda_1^{\ast},\ldots,\lambda_4^{\ast}\}$
of $H^1(X,\Z)$), then by the description of the Hodge structure on
the first cohomology of an abelian variety given in the previous
section,
$\phi_i$ is the VHS associated to the family
$p_i:\X_D\to\Delta^{\ast}(M)$.

\subsection{The limit mixed Hodge structure}\label{subsection35}

The family $p_i:\X_D\to\Delta^{\ast}(M)$ ($i=1$,2) induces a MHS on
$H^1(X_t)$ ($t\in\Delta^{\ast}(M)$ arbitrary but fixed; see
\cite[Chapter IV]{G}). Its weight filtration is of the form
$0=W_{-1}\subset W_0 \subset W_1 \subset W_2=H^1(X_t)$. It is is
given by
$$
     W_0 = \im N\text{ and }
     W_1 = \Ker N,
$$ where $N=\log T=T-\mbox{\boldmath$1$}=\left(\smallmatrix 0 &
\eta\\  0 & 0\endsmallmatrix\right)$. (We have written $T$ instead of
$T_i$ etc.) The limit Hodge filtration $F_{\infty}^{\bullet}\subset
H^1(X_t)$ is calculated as follows: first define $\tilde{\psi}:
\h_1(M) \to \check{\D}_2\subset\Grass(2,4)$ (where $\check{\D}_2$ is
the compact dual of $\D_2$) by
$$
\sideset{^t}{}{\tilde{\psi}}(\tau)
=T^{-\tau}\sideset{^t}{}{\tilde{\phi}}(\tau)=
(\mbox{\boldmath$1$}-\tau N)
\sideset{^t}{}{\tilde{\phi}}(\tau).
$$ By construction $\tilde{\psi}(\tau+1)=\tilde{\psi}(\tau)$, hence
$\tilde{\psi}$ descends to a map $\psi: \Delta^\ast(M)\to
\check{\D}_2.$  By the work of Griffiths and Schmid $\psi$ extends to
$\Delta$ and $F_{\infty}^{\bullet}:=\psi(0)$ together with the
monodromy weight filtration defines a MHS on $H^1(X_t)$.

{}From now on  we will write $e_i$ instead of $\lambda_i^{\ast}$ for
the basis vectors of $H^1(X,\Z)$. In our case, an immediate
calculation shows that $\tilde{\psi}$ is constant  and
$$
    W_0=[e_1],\quad
    W_1=[e_1, e_2, e_4],\quad
    F_{\infty}^1=\sqper{0}{\tau_3}.
$$  for the family $p_1:\X_D\to\h_1(M)$ and
$$
    W_0=[e_2],\quad
    W_1=[e_1, e_2, e_3],\quad
    F_{\infty}^1=\sqper{\tau_1}{0}.
$$  for the family $p_2:\X_D\to\h_1(M)$.

We shall now restrict our attention to $W_1$. Here the Hodge
filtration is given by
$$
     F_{\infty}^1\cap W_1 =
     [\tau_2 e_1 + \tau_3 e_2 + pe_4]
$$ resp.
$$
     F_{\infty}^1\cap W_1 =
     [\tau_1 e_1 + \tau_2 e_2 + e_3].
$$ By abuse of notation we will denote this by $F_{\infty}^1$, too.
Finally, recall that all elements of the nilpotent orbit
$\{T^{\tau}F_{\infty}^1\}$ of Hodge filtrations on $W_2=H^1(X_t)$
induce the same MHS on $W_1$ (cf. \cite[p.84]{G}).

\subsection{Calculation of $\mbox{Gr}_1^W$ and the extension
class}\label{subsection36}

In the previous section we computed the limit MHS on $H^1(X_t)$
determined by $p_i:X\to\Delta^*(M)$ ($i=1,2$). We now want to compute
$\operatorname{Gr}_1^W$ and the extension class.

For any integer $n$ and $\tau\in\h_1$ we define
$$ X_{\tau,n}=V/\Lambda=\C f/(\Z\mu_1+\Z\mu_2)
$$ where $\mu_1=\tau f$ and $\mu_2=nf$.

\begin{pro}\label{P1}
$\operatorname{(i)}$ In the central case ($i=1$)
$\operatorname{Gr}_1^W=H^1(C)$ where
$C=X_{\tau_3,p}$.

\noindent $\operatorname{(ii)}$ In the peripheral case ($i=2$)
$\operatorname{Gr}_1^W=H^1(C)$ where
$C=X_{\tau_1,1}$.
\end{pro}
\begin{pf} (i) Setting $\mu_1^*=e_2$ and $\mu_2^*=e_4$ we have
$$ F^1=[f^*]=[\tau_3\mu_1^*+p\mu_2^*]=[\tau_3e_2+pe_4]
$$ and the claim follows immediately from the  calculations of
\ref{subsection35}.

(ii) Setting $\mu_1^*=e_1$ and $\mu_2^*=e_3$ we have
$$ F^1=[f^*]=[\tau_1\mu_1^*+\mu_2^*]=[\tau_1e_1+e_3]
$$ and the result follows as before.
\end{pf}

Next we want to calculate the class of the extensions
\begin{equation}\label{2} 0\longrightarrow
W_0\stackrel{i}{\longrightarrow}W_1\stackrel{\pi}{\longrightarrow}
\operatorname{Gr}_1^W\longrightarrow 0.
\end{equation}

The recipe for finding the extension class of (\ref{2}) demands that
we choose an integral retraction $r: W_1\to W_0$ and then transfer
$F_{\infty}^1$ to $W_0\oplus \Gr_1^W$ via $(r,\pi)$. In the central
case we take
$r(e_1)=e_1$ and $r(e_2)=r(e_4)=0$ and in the peripheral case
$r(e_2)=e_2$ and $r(e_1)=r(e_3)=0$. Then, in both cases,
$(r,\pi)=\id_{\Gr_1^W}$, hence $\tilde{F}_{\infty}^1=F_{\infty}^1$.
Next we define the trivial extension
$$
     \tilde{F}_0^1 = F^1(W_0)\oplus F^1(\Gr_1^W)
          = F^1(\Gr_1^W) = [\tau_3 e_2 +pe_4]
$$ for the first family and similarly
\[
   \tilde{F}_0^1 =[\tau_1 e_1 + e_3]
\] for the second family. The extension class $e\in
J^0\Hom(\Gr_1^W,W_0)$ is represented by any map
$\psi: \Gr_1^W \to W_0$ such that
$$
     \mwp{W_0}{\Gr_1^W} \tilde{F}_0^1 =\tilde{F}_{\infty}^1.
$$ For the first family we can take $\psi=(\tau_2/p)
e_4^{\ast}\otimes e_1$ and for the second family $\psi=\tau_2
e_3^{\ast} \otimes e_2$.

After identifying $W_0=\Z$ and $\Gr_1^W=H^1(C)$, where $C$ is as in
proposition~\ref{P1}, the extension class $e=[\psi]$ is well defined
in
\begin{align*}
\operatorname{Ext}(\operatorname{Gr}_1^W,W_0)
          &= J^0 \Hom(\Gr_1^W,W_0)\\
          &= H^1(C)^\ast / (H^{0,1}(C))^\ast + H^1(C,\Z)^\ast\\
          &= H^{1,0}(C)^\ast / H_1(C,\Z)\\
          &=\Alb(C)
\end{align*} and is represented by $\psi|_{H^{1,0}}$.

\begin{pro}\label{P2} Under the identification
$\operatorname{Ext}(\operatorname{Gr}_1^W,W_0)=\operatorname{Alb}(C)$
the extension class of (\ref{2}) corresponds to
$[\tau_2f]\in\operatorname{Alb}(C)$.
\end{pro}
\begin{pf} Notice that if
$C=V/\Lambda$, then $H^1(C,\C)=\Hom_{\R}(V,\C)$ and
\[
\begin{CD}
     H^1(C,\C) @>{\int_{0}^{[v]}}>> \C/{\operatorname{periods}}\\
     @| @AAA\\
     \Hom_{\R}(V,\C) @>{\operatorname{ev}(v)}>> \C
\end{CD}
\] commutes.

The extension class $e\in\Ext(\Gr_1^W,W_0)$  corresponds to $[\tau_2
f]\in\Alb(C)$ since
\begin{align*}
     \int_0^{[\tau_2 f]} (\tau_3 e_2 + p e_4)
          &= \int_0^{[\tau_2 f]} (f^{\ast})\\
          &=\tau_2\\
          &=\psi(\tau_3 e_2 + p e_4)
\end{align*} in the central case and
\begin{align*}
     \int_0^{[\tau_2 f]} (\tau_1 e_1 + e_3)
          &= \int_0^{[\tau_2 f]} (f^{\ast})\\
          &=\tau_2\\
          &=\psi(\tau_1 e_2 + e_3)
\end{align*} in the peripheral case.
\end{pf}

\subsection{The number of components}\label{sec47}

In this section we will proof the following proposition.

\begin{pro}\label{tnoc} Let $\pi\colon  X \to \Delta$ be a semi-stable
degeneration of $(1,p)$-polarized abelian surfaces with (canonical)
level structure. Assume that the central fibre $X_0$ is a cycle of
$N$ ruled surfaces. Then the number $N$ is determined by the
polarized VHS (PVHS) on $\Delta^\ast$.
\end{pro}

The polarization on the VHS is induced by the polarization $\cL$ in
the following way: For $Y\subset X$ define
\[
   c_Y = c_1(\cL|_Y) \in H^2(Y)
\] and
\[
     Q_Y \colon  (a,b) \mapsto a \cup b \cup c_Y\colon
          H^1(Y) \times H^1(Y) \to H^4(Y).
\] We are mostly interested in the case where $Y$ is $X$,
$X_0$ or $X_t$. We write $c_0$ for $c_{X_0}$ and $Q_0$ for $Q_{X_0}$.
Similarly, we write $c_t$ for $c_{X_t}$ and $Q_t$ for $Q_{X_t}$.

The existence of a global level structure is needed to insure that
$\deg(\cL|_{F_k})$, where $F_k$ is the fibre of $Y_k$, is independent
of $k$. Indeed, by lemma 4 of \S2.2 it is always 1.

For the rest of the section $\pi\colon  X \to \Delta$ is as in
proposition \ref{tnoc}. As in \S3.1 we identify $C$ with $C_0=Y_0\cap
Y_1$. Let $Y_k$ and $\tau_k\colon  Y_k \to C$ be as in \S3.1 and let
$i\colon  C\into X_0$ and $j_k\colon Y_k \into X_0$ be the inclusion
maps.

\begin{lem}\label{cpt} Under the isomorphism
$(\tau_k^\ast)\colon  H^1(C) \to \Gr_1^W$, the induced map
$(j_k^\ast)\colon  W_1 \to \Gr_1^W$ corresponds to the map
$i^\ast\colon H^1(X_0)\to H^1(C)$ induced by the inclusion
$i\colon  C\into X_0$.
\end{lem}

\begin{pf} This follows immediately from
$j_k^\ast=\tau_k^\ast i^\ast\colon  H^1(X_0) \to H^1(Y_k)$, which in
turn follows from the following facts:
$\rho_k^\ast$ is the inverse of $\alpha_k^\ast$,
$\sigma_k^\ast$ is the inverse of $\beta_k^\ast$,
$i=i_0=j_1\beta_0$, and $j_{k+1}\beta_k= j_k\alpha_k$.
\end{pf}

\begin{pf1} The PVHS on $\Delta^\ast$ determines the bilinear form
$Q_t$ on $W_2=H^1(X_t)$. Since the maps on cohomology induced by an
inclusion $Y'\subset Y$ are compatible with $Q_Y$ and $Q_{Y'}$, $Q_t$
induces
$Q_0$ on $H^1(X_0)=W_1\subset W_2$, i.e., we have a commutative
diagram
$$
\begin{CD}
     H^1(X_t) \otimes H^1(X_t) @>{Q_t}>> H^4(X_t)\\
     @A{k_t^\ast\times k_t^\ast}AA @AA{k_t^\ast}A\\
     H^1(X) \otimes H^1(X) @>{Q}>> H^4(X)\\
     @V{k_0^\ast\times k_0^\ast}VV @VV{k_0^\ast}V\\
     H^1(X_0) \otimes H^1(X_0) @>{Q_0}>> H^4(X_0).
\end{CD}
$$ For $t\in\Delta$ let $k_t\colon X_t\into X$ be the inclusion map
and let
$$
     \delta\colon H^4(X_0) \to H^4(X_t)
$$ be $k_t^\ast (k_0^\ast)^{-1}$ (recall that $k_0^\ast$ is an
isomorphism). Since the monodromy is compatible with
$Q_t$, i.e.,
$
     Q_t(Ta,Tb)= Q_t(a,b),
$ so is $N=T-\id$, i.e.
$
     Q_t(Na,b)=Q_t(a,Nb).
$ Since
     $W_1=$ Ker $N$ and $W_0=$ Im $N$
$
     Q_t(W_1,W_0)=0.
$ In other words, there exists a bilinear form $q$ on $H^1(C)$ making
the diagram
\begin{center}
\unitlength1pt
\begin{picture}(190,70)
\put(0,0){\makebox(90,10){$H^1(C)\times H^1(C)$}}
\put(0,50){\makebox(90,10){$H^1(X_0)\times H^1(X_0)$}}
\put(130,50){\makebox(60,10){$H^4(X_t)=\Z$}}
\put(45,40){\vector(0,-1){20}}
\put(5,20){\makebox(40,20){$i^\ast\times i^\ast$}}
\put(100,55){\vector(1,0){20}}
\put(75,20){\vector(2,1){40}}
\put(100,18){$q$}
\put(100,55){\makebox(20,20){$\delta Q_0$}}
\end{picture}
\end{center} commute. Furthermore, since $i^\ast\colon  H^1(X_0) \to
H^1(C)$ is surjective,
$q$ is determined by $\epsilon Q_0$ where
$$
\varepsilon: H^2(C)\longrightarrow H^4(X_t)
$$ is the isomorphism which sends the generator of $H^2(C)$
corresponding to the canonical orientation of $C$ to the same thing
on $X_t$. Let
$$
\gamma:H^1(X)\times H^1(C)\longrightarrow H^2(C)
$$ be the intersection form. We will show in corollary \ref{mcc} that
$q=N\varepsilon\gamma$. This means that the data over
$\Delta^*$ determine $N$.
\end{pf1}

\medskip

For convenience, we shall first assume that $N\ge2$. The case $N=1$
requires minor modifications and is dealt with later.

Since the line bundle determining $Y_k$ has degree $0$, there exists
a continuous map $q_k\colon  Y_k \to \PR^1$ such that
$$
     (\tau_k,q_k)\colon  Y_k \to C \times \PR^1
$$ is a homeomorphism. Let $f\in H^2(\PR^1)$ be the canonical
generator, determined by the complex structure. We define $\beta$ by
demanding that the diagram
$$
\begin{CD} H^2(C) @>{\beta}>> H^4(X_0)\\ @V{(\id,\cdots,\id)\otimes
f}VV @VV{(j_k^\ast)}V\\
\oplus H^2(C)\otimes H^2(\PR^1) @>>{(\tau_k^\ast \cup q_k^\ast)}>
\oplus H^4(Y_k)
\end{CD}
$$ be commutative. (Notice that $(j_k^\ast)$ is an isomorphism.)

\begin{lem}\label{mcl1} The diagram
$$
\begin{CD}
     H^1(X_0)\times H^1(X_0) @>{Q_0}>> H^4(X_0)\\
     @V{i^\ast\times i^\ast}VV @AA{\beta}A\\
     H^1(C) \times H^1(C) @>>{\gamma}> H^2(C)
\end{CD}
$$ commutes.
\end{lem}

\begin{pf} By definition of $\beta$ we have to show that
$$
\begin{CD}
     H^1(X_0)\times H^1(X_0) @>{Q_0}>> H^4(X_0) @>{(j_k^\ast)}>>
          \oplus H^4(Y_k)\\
     @V{i^\ast\times i^\ast}VV @. @AA{(\tau_k^\ast\cup q_k^\ast)}A\\
     H^1(C) \times H^1(C) @>>{(\gamma,\ldots,\gamma)}>
          \oplus H^2(C) @>>{{}\otimes f}>
     \oplus H^2(C) \otimes H^2(\PR^1)
\end{CD}
$$ commutes. In lemma 4 of \S2 we showed that
$\deg(\cL|_{F_k})=1$ for all $k$, where $F_k$ is the ruling of $Y_k$.
Hence we can write
$
     j_k^\ast (c_0) = \tau_k^\ast (d_k) + q_k^\ast (f)
$ for some $d_k\in H^2(C)$. Since $j_k^\ast=\tau_k^\ast i^\ast$ on
$H^1(X_0)$ (see proof of lemma \ref{cpt})
$$
     j_k^\ast  Q_0(a,b) = \tau_k^\ast i^\ast (a\cup b) \cup
          (\tau_k^\ast (d_k) +  q_k^\ast (f))
          = \tau_k^\ast i^\ast (a \cup b) \cup q_k^\ast (f),
$$ because $i^\ast(a\cup b)\cup d_k\in H^4(C)=0$.
\end{pf}

The orientations of the components of the singular fibre are
compatible with the orientation of the general fibre in the following
sense:

\begin{lem}\label{mcl2} If $N$ is the number of components of $X_0$,
then $\delta\beta=N\epsilon$.
\end{lem}

\begin{pf} By definition of $\beta$, it suffices to show that the
diagram
$$
\begin{CD}
     H^4(X_0) @<{k_0^\ast}<{\sim}< H^4(X)\\
     @V{(j_k^\ast)}V{\sim}V @VV{k_t^\ast}V\\
     \oplus H^4(Y_k) @>>{(\id,\ldots,\id)}> H^4(X_t)
\end{CD}
$$ commutes. The inverse of $k_0^\ast$ is induced by a retraction
$r\colon X \to X_0$, which exhibits $X_0$ as a strong deformation
retract of $X$ \cite{Cl, P, St}. This retraction restricts to
$r_t\colon X_t\to X_0$. Let $Z_k\subset X_t$ be the inverse image of
$Y_k$. Then
$S_k:=Z_k\cap Z_{k+1}$ is the inverse image of
$T_k:=Y_k\cap Y_{k+1}$. Pick any ${\kappa}\in\{1,\ldots,N\}$. Let
$Y'=Y_{{\kappa}}$, $Y''=\cup_{k\neq {\kappa}} Y_k$,
$Z'=Z_{{\kappa}}$ and $Z''=\cup_{k\neq {\kappa}} Z_k$. Let $S=Z'\cap
Z''=S_{{\kappa}-1}\cup S_{{\kappa}}$ and
$T=Y'\cap Y''=T_{{\kappa}-1}\cup T_{{\kappa}}$.

We will need below that $H^4(Z')=H^4(Z'')=0$. To see this, we use that
the retraction $r_t\colon Z_k\to Y_k$ is the real oriented blow up
along
$T_{k-1}\cup T_k$ \cite[p.~36]{P}. The \lq\lq exceptional
divisors\rq\rq\ $S_k=r_t^{-1}(T_k)$ are
$S^1$-bundles over $T_k$. They are trivial, because the triple
$(Y_k, T_{k-1}, T_k)$ is homeomorphic to
$(B\times S^2, B\times\{0\}, B\times\{\infty\})$, where $B=S^1\times
S^1$. It follows that the triple
$(Z_k,S_{k-1},S_k)$ is homeomorphic to
$(B\times S^1\times [0,1], B\times S^1\times \{0\}, B\times S^1\times
\{1\})$. In particular, $H^4(Z_k)=0$ and
$H^3(Z_k)\to H^3(S_k)$ and
$H^3(Z_k)\to H^3(S_{k-1})$ are surjective for all $k$. Mayer-Vietoris
now implies that
$H^4(Z'')=\oplus_{k\neq {\kappa}} H^4(Z_k)$ and this proves the claim.

Consider the commutative diagram
$$
\begin{CD}
     H^4(Y') @<<< H^4(X_0) @>{r_t^\ast}>> H^4(X_t)\\
     @A{\sim}AA @AAA @AA{\sim}A\\
     H^4(Y',T) @<{\sim}<< H^4(X_0,Y'') @>{\sim}>{r_t^\ast}>
H^4(X_t,Z'')
\end{CD}
$$ The map $H^4(X_0,Y'')\to H^4(X_t,Z'')$ is an isomorphism by
Alexander duality since $r_t$ is an homeomorphism
$X_0\setminus Y''\to X_t\setminus Z''$ (cf.\ \cite[p.23]{La}).
Similarly, $ H^4(X_0,Y'')\to H^4(Y',T)$ is an isomorphism. The left
vertical map is an isomorphism because $T$ has topological dimension
2. The right vertical map is surjective because
$H^4(Z'')=0$, as we have seen above. But a surjective map form $\Z$
to $\Z$ is automatically injective.

It follows from \cite[p.99]{Gre} that the map
$H^4(Y')\to H^4(X_0)$ via $H^4(Y',T)$ and $H^4(X_0,Y'')$ composed
with the Mayer-Vietoris isomorphism
$H^4(Y')\oplus H^4(Y'') \cong H^4(X_0)$ is just the inclusion of
$H^4(Y')$ into $H^4(Y')\oplus H^4(Y'')$. We have to show that the
composition of this map with $r_t^\ast\colon H^4(X_0)\to H^4(X_t)$
maps the orientation class to the orientation class. As before, we
may identify the relative cohomology groups with the cohomology with
compact supports  of the respective complements. But $r_t\colon
X_t\setminus Z''\to Y'\setminus T$ is a complex isomorphism. In
particular, it respects the orientation.
\end{pf}

\begin{cor}\label{mcc} $q=N\epsilon\gamma$.
\end{cor}

\begin{pf} By lemma \ref{mcl1}
$q=\delta\beta\gamma$. Now the result follows immediately from the
previous lemma \ref{mcl2}.
\end{pf}

Finally, we consider the case $N=1$. We proceed as in the case
$N\ge 2$. First we have to define
$\beta\colon H^2(C)\to H^4(X_0)$. To do this, notice that we have an
exact Mayer-Vietoris sequence
$$
      H^{q-1}(\tilde{C})\to
     H^q(X_0)\to H^q(\tilde{X_0})\oplus H^q(C)\to H^q(\tilde{C}),
$$ where $\tilde{X_0}$ is the normalization of $X_0$ and
$\tilde{C}=\pi^{-1}(C)$ (see \cite[p.120]{C}). Hence the
normalization map
$\pi\colon \tilde{X_0}\to X_0$ induces an isomorphism
$H^4(X_0)\to H^4(\tilde{X_0})$. We identify $H^4(X_0)$ with
$H^4(\tilde{X_0})$ via $\pi^\ast$ and define $\beta$ by demanding
that it maps the orientation class of $C$ to the orientation class of
$\tilde{X_0}$. With this definition of $\beta$ and using that the
degree $\deg(\pi^\ast(\cL)|_F)$ of
$\cL$ along the fibre $F$ in the normalization is 1 again by lemma 4
of \S2, the diagram in lemma~\ref{mcl1} again commutes. The proof is
an easy adaption of the proof given for $N\ge2$ and is left to the
reader.

To prove lemma \ref{mcl2} for $N=1$, one argues exactly as in the
case $N\ge2$ using the commutative diagram
$$
\begin{CD}
     H^4(\tilde{X_0}) @<<< H^4(X_0) @>{r_t^\ast}>> H^4(X_t)\\
     @AAA @AAA @AAA\\
     H^4(\tilde{X_0},\tilde{C}) @<<< H^4(X_0,C)
          @>{r_t^\ast}>> H^4(X_t,S),
\end{CD}
$$ where $S=r_t^{-1}(C)$. (This time, all maps are isomorphisms.)

\section{A uniqueness result}\label{section5}

Here we combine the results of sections \ref{section3} and
\ref{section4} to prove a uniqueness result for degenerate abelian
surfaces.

In (\ref{subsection34}) we discussed central and peripheral corank 1
boundary points of ${\cal A}^*(1,p)$. Note that the universal family
$\pi:X_D\to\h_2$ descends to a family $\pi':X_D\to e(\h_2)$. By
\cite[part~II.4]{HKW} this family can be extended to a family
$\bar\pi:\bar{X}_D\to X(l)$. The fibres of $\bar{X}_D$ over the
''boundary''
$\{0\}\times\C\times\h_1$, resp.
$\h_1\times\C\times\{0\}$ are cycles of elliptic ruled surfaces with
$N=1$, resp. $N=p$.

\begin{dfn}\label{D} Let $[q]\in{\cal A}^*(1,p)$ be a corank 1
boundary point. We say that a surface $X_0$ is a {\em degenerate
abelian surface associated to} $[q]$ if the following holds: There
exists a 1-parameter degeneration of polarized abelian surfaces with
(canonical) level structure
$p:X\to\Delta$  as in section 2.2 with central fibre $p^{-1}(0)=X_0$
and an  embedding
$f:\Delta\to X(l)$ such that:

\noindent (i) $f(\Delta)$ meets $\{0\}\times\C\times\h_1$ resp.
$\h_1\times\C\times\{0\}$ transversely in the point $q=f(0)$ which is
mapped to $[q]$.

\noindent (ii) The restriction of $p:X\to\Delta$ to $\Delta^*$ is the
pull back of the family
$\pi':X_D\to e(\h_2)$ via $f$.
\end{dfn}

\begin{rem}\normalshape The reason why we work with $X(l)$ rather
than with the compactification ${\cal A}^*(1,p)$ itself is that the
map $\h_2\to\cal A(1,p)$ is branched  over two Humbert surfaces $H_1$
and $H_2$. Near these surfaces, resp. their closure in ${\cal
A}^*(1,p)$ one has no universal family.
\end{rem}

\begin{lem}\label{lem25} If $X_0$ is a degeneration with only local
normal crossing singularities associated  to a corank 1 boundary
point $[q]$, then  the MHS on $H^1(X_0)$ is determined by this point.
\end{lem}

\begin{pf} This follows from the Local Invariant Cycle theorem
\cite[Ch.~VI]{G} (which also holds when we have local normal crossing
rather than global normal crossing) together with the existence of the
extended family
$\bar{\pi}:\bar{X}_D\to X(l)$: Indeed, compare  the family
$p:X\to\Delta$ as in definition (\ref{D}) with the family
$f^{\ast}\bar{X}_D\to\Delta$. On $\Delta^{\ast}$ these two families
agree (they are both the pull back of the universal family), hence
they determine the same VHS and thus the same limit MHS  on
$H^1(X_t)$ ($t\in\Delta^{\ast}$ arbitrary but fixed). But by the
Local Invariant Cycle theorem the MHS on the central fibre is
determined by this limit MHS at least over $\Q$: more precisely,
$H^1(X_0)=W_1(H^1(X_t))=H^1(\bar{\pi}^{-1}(q))$ as $\Q-$MHS. It
remains to show that
$$ 0\longrightarrow H^1(X)\longrightarrow
H^1(X_t)\stackrel{N}{\longrightarrow} H^1(X_t)
$$ is exact over $\Z$. For $n\gg0$ the zero locus $Y\subset X$ of a
general element of
$\Gamma (X,\cL^n)$ is a semi-stable degeneration of the smooth curve
$Y_t=Y\cap X_t$. By the Picard-Lefschetz theorem \cite[theorem
III.~14.1]{BPV}
$$ 0\longrightarrow H^1(Y)\longrightarrow
H^1(Y_t)\stackrel{N}{\longrightarrow} H^1(Y_t)
$$ is exact over $\Z$. By the Lefschetz hyperplane theorem
$H^1(X_t)\longrightarrow H^1(Y_t)$ is injective, hence so is
$H^1(X)\longrightarrow H^1(Y)$. We have now reduced our problem to
showing that $H^1(Y)/H^1(X)$ is torsion free. Since
$H^1(Y)/H^1(X)\hookrightarrow H^2(X,Y)$ it is sufficient to show that
$H^2(X,Y)$ or, equivalently $H_1(X,Y)=\text{Coker }
(H_1(Y)\longrightarrow H_1(X))$ is torsion free. Since dim ker $N=3$
the central fibre $X_0$ is not smooth abelian. Hence by Persson's
theorem
$X_0$ is a cycle of elliptic ruled surfaces.  Assume that $X_0$ has
at least two components. (The remaining case is analogous.) Let $C$
be one of the double curves and
$X_i$ one of the components of $X_0$. Then
$$ H_1(X)=H_1(X_0)\cong H_1(C)\oplus \Z c
$$ where $c\in H_1(X_0)$ is as in 3.3. Let $Y_i=X_i\cap Y$. Then
$X_i$ and
$Y_i$ are smooth and irreducible. We have a commutative diagram
$$
\begin{CD}
     H_1(C)= @. H_1(X_i) @>>> H_1(X)\\
                       @. @AAA @AAA\\
                       @. H_1(Y_i)  @>>> H_1(Y).
\end{CD}
$$ Since $c$ is clearly induced from $Y$ up to $H_1(C)$, it is
sufficient to prove that
$$ H_1(Y_i)\longrightarrow H_1(X_i)
$$ is surjective. But this follows from the Lefschetz hyperplane
theorem since
$Y_i$ is ample on $X_i$. This shows that $H^1(X_0)=W_1(H^1(X_t))$
depends only on the point$[q]$.
\end{pf}

\begin{thm}\label{T} Let
$X_0$ be a degeneration associated to a corank 1 boundary point $[q]$
of
${\cal A}^*(1,p)$ with only local normal crossing singularities and
no triple points. Then
$X_0$ is a cycle of ruled surfaces (possibly with only one component)
over an elliptic curve $C$. There exists a $\cL\in\Pic^0(C)$ such that
all components are isomorphic to
$\PR^1(\cO\oplus\cL)$. If $q$ is a central boundary point then $X_0$
is completely determined by $q$. If $q$ is peripheral, then the
number $N$ of components, the base curve $C$ andthe shift $s$ are
uniquely determined by the point $[q]$. The line bundle
$\cL\in\Pic^0(C)$ determining the elliptic ruled surfaces of the
cycle is determined (at least) up to a $p$-torsion point of
$\Pic^0(C)$. More precisely:
\begin{enumerate}
\item If $[q]=[\tau_2,\tau_3]$ is a central boundary point, then
$N=1$,
$C=X_{\tau_3,p}$, $s=[\tau_2]$ and $\cL=\cO(s-O)$.
\item If $[q]=[\tau_1,\tau_2]$ is a peripheral boundary point, then
$N=p$,
$C=X_{\tau_1,1}$, $s=[\tau_2]$ and there exists a point $s'\in C$
such that $\cL=\cO(s'-O)$.
\end{enumerate}
\end{thm}

\begin{pf} We first note that the pair $(X_{\tau_3,p},[\tau_2])$,
resp.
$(X_{\tau_1,1},[\tau_2])$ only depends on the point $[q]$ but not on
its representative \cite[part~I.3]{HKW}.

If $X_0$ has only global normal crossings then Persson's result
\cite[proposition 3.3.1]{P} says that $X_0$ is a cycle of elliptic
ruled surfaces.

The previous lemma asserts that the MHS on $H^1(X_0)$ is determined
by $q$. Furthermore,  this MHS is the weight 1 part of the limit MHS
on  $H^1(X_t)$ of the family $p_1$ (resp.
$p_2$) defined in \S4.4 in the central (resp. peripheral) case.
Indeed, $p_1$ (resp. $p_2$) is the restriction of the universal family
to a small disk $\Delta\times (\tau_2,\tau_3)$ (resp.
$(\tau_1,\tau_2)\times\Delta$). By lemma~\ref{wf} the base curve $C$
is determined by the equation $H^1(C)=
\Gr_1^W(H^1(X_0))=\Gr_1^W(H^1(X_t))$ and this equation was solved in
proposition~\ref{P1}. Corollary~\ref{cor17} shows how to recover the
shift $s$ from
$H^1(X_0)=W_1(H^1(X_t))$ and the result for the families $p_1$ and
$p_2$ is given in proposition~\ref{P2}.

Applying proposition \ref{tnoc} of \S\ref{sec47} to the family $X$
defining
$X_0$ and to the family constructed in \cite{HKW}, one sees that the
number of components of $X_0$ is independent of the choice of the
disk $\Delta\into X(l)$ (cf.\ proof of lemma \ref{lem25}). But it
follows from our calculations in \S4 that the polarization induces
$\left(\begin{array}{cc} 0& -1\\1&0\end{array}\right)$ (with respect
to the basis
$\{e_2,e_4\}\bmod e_1$ of \S\ref{subsection35}) in the central case
and
$\left(\begin{array}{cc} 0& -p\\p&0\end{array}\right)$ (with respect
to the basis $\{e_1,e_3\}\bmod e_2$) in the peripheral case.  It
follows from proposition \ref{tnoc} of \S\ref{sec47} that the number
of components is 1 and $p$ respectively.

The statement about the form of $\cL\in\Pic^0(C)$ now follows directly
from proposition 6 of \S2.2
\end{pf}

\begin{rem}\normalshape Note that in the proof of theorem \ref{T}
the level structure is not needed to determine $C$ and $s$.
\end{rem}

\begin{rem} \normalshape It seems reasonable to expect that also in
the peripheral case
$\cL$ is completely determined by the boundary point $q$. To show
this in the spirit of this paper, one probably has to translate the
notion of (canonical) level structure to Hodge structures.
\end{rem}

\begin{rem}\normalshape This result is compatible with \cite{HKW}.
Notice that we only used the existence of the degeneration
$\bar{X}_D$ as constructed there in order to prove our uniqueness
result, but not the precise description of its singular fibres.
\end{rem}

\vspace{2cm}

\parbox{5cm}{K.~Hulek\\ Institut f\"ur Mathematik\\ Universit\"at
Hannover\\ Postfach 6009\\ 30060 Hannover\\ Germany}\hspace{1cm}
\parbox{5cm}{J.~Spandaw\\ Mathematisches Institut\\ Universit\"at
Bayreuth\\ Postfach 101251\\ 95440 Bayreuth\\ Germany}
\end{document}